# Multivariate volatility models

## Ruey S. Tsay[1]

*University of Chicago*

**Abstract:** Correlations between asset returns are important in many financial applications. In recent years, multivariate volatility models have been used to describe the time-varying feature of the correlations. However, the curse of dimensionality quickly becomes an issue as the number of correlations is $k(k-1)/2$ for $k$ assets. In this paper, we review some of the commonly used models for multivariate volatility and propose a simple approach that is parsimonious and satisfies the positive definite constraints of the time-varying correlation matrix. Real examples are used to demonstrate the proposed model.

## 1. Introduction

Let $r_t = (r_{1t}, \ldots, r_{kt})'$ be a vector of returns (or log returns) of $k$ assets at time index $t$. Let $F_{t-1}$ be the sigma field generated by the past information at time index $t-1$. We partition the return $r_t$ as

$$r_t = \mu_t + e_t, \tag{1}$$

where $\mu_t = E(r_t|F_{t-1})$ is the conditional mean of the return given $F_{t-1}$ and $e_t$ is the innovation (or noise term) satisfying $e_t = \Sigma_t^{1/2}\epsilon_t$ such that

$$\text{Cov}(e_t|F_{t-1}) = \text{Cov}(r_t|F_{t-1}) = \Sigma_t, \tag{2}$$

where $\epsilon_t = (\epsilon_{1t}, \ldots, \epsilon_{kt})'$ is a sequence of independently and identically distributed random vectors with mean zero and identity covariance matrix, and $\Sigma_t^{1/2}$ is the symmetric square-root matrix of a positive-definite covariance matrix $\Sigma_t$, that is, $\Sigma_t^{1/2}\Sigma_t^{1/2} = \Sigma_t$. In the literature, $\Sigma_t$ is often referred to as the volatility matrix. Volatility modeling is concerned with studying the evolution of the volatility matrix over time. For asset returns, behavior of the conditional mean $\mu_t$ is relatively simple. In most cases, $\mu_t$ is simply a constant. In some cases, it may assume a simple vector autoregressive model. The volatility matrix $\Sigma_t$, on the other hand, is much harder to model, and most GARCH studies in the literature focus entirely on modeling $\Sigma_t$.

The conditional covariance matrix $\Sigma_t$ can be written as

$$\Sigma_t = D_t R_t D_t \tag{3}$$

where $D_t$ is a diagonal matrix consisting of the conditional standard deviations of the returns, i.e., $D_t = \text{diag}\{\sqrt{\sigma_{11,t}}, \ldots, \sqrt{\sigma_{kk,t}}\}$ with $\sigma_{ij,t}$ being the $(i,j)$th element of $\Sigma_t$, and $R_t$ is the correlation matrix.

---

In recent years, many studies extend the univariate generalized autoregressive conditional heteroscedastic (GARCH) model of Bollerslev [2] to the multivariate case for modeling the volatility of multiple asset returns; see the recent article [1] for a survey. Multivariate volatility models have many important applications in finance and statistics. They can be used to study the correlations between asset returns. These correlations play an important role in asset allocation, risk management, and portfolio selection. There are two major difficulties facing the generalization, however. First of all, the dimension of volatility matrix increases rapidly as the number of asset increases. Indeed, there are $k(k+1)/2$ variances and covariances for $k$ asset returns. Second, for asset returns the covariance matrix is time-varying and positive definite. Many of the multivariate volatility models proposed in the literature fail to satisfy the positive-definite constraints, e.g., the diagonal VEC model [3], even though they are easy to understand and apply.

The goal of this paper is to propose a simple approach to modeling multivariate volatility. The proposed model is kept parsimonious in parameterization to overcome the difficulty of curse of dimensionality. In addition, a simple structure equation is imposed to ensure that the resulting time-varying covariance matrices are positive definite. On the other hand, the proposed model is not very flexible and may encounter lack of fit when the dimension is high. To safe guard against model inadequacy, we consider model checking using some bootstrap methods to generate finite-sample critical values of the test statistics used.

The paper is organized as follows. In Section 2, we briefly review some of the multivariate volatility models relevant to the proposed model. Section 3 considers the proposed model whereas Section 4 contains applications to daily returns of foreign exchange rates and U.S. stocks. Section 5 concludes.

## 2. A brief review of vector volatility models

Many multivariate volatility models are available in the literature. In this section, we briefly review some of those models that are relevant to the proposed model. We shall focus on the simple models of order $(1,1)$ in our discussion because such models are typically used in applications and the generalization to higher-order models is straightforward. In what follows, let $a_{ij}$ denote the $(i,j)$th element of the matrix $A$ and $u_{it}$ be the $i$th element of the vector $u_t$.

*VEC model.* For a symmetric $n \times n$ matrix $A$, let vech($A$) be the half-stacking vector of $A$, that is, vech($A$) is a $n(n+1)/2 \times 1$ vector obtained by stacking the lower triangular portion of the matrix $A$. Let $h_t = \text{vech}(\Sigma_t)$ and $\eta_t = \text{vech}(e_t e_t')$. Using the idea of exponential smoothing, Bollerslev et al. [3] propose the VEC model

$$(4) \qquad h_t = c + A\eta_{t-1} + Bh_{t-1}$$

where $c$ is a $k(k+1)/2$-dimensional vector, and $A$ and $B$ are $k(k+1)/2 \times k(k+1)/2$ matrices. This model contains several weaknesses. First, the model contains $k(k+1)[k(k+1)+1]/2$ parameters, which is large even for a small $k$. For instance, if $k=3$, then the model contains 78 parameters, making it hard to apply in practice. To overcome this difficulty, Bollerslev et al. [3] further suggest that both $A$ and $B$ matrices of Eq. (4) are constrained to be diagonal. The second weakness of the model is that the resulting volatility matrix $\Sigma_t$ may not be positive definite.

*BEKK model.* A simple BEKK model of Engle and Kroner [5] assumes the form

$$\Sigma_t = C'C + A'e_{t-1}e'_{t-1}A + B'\Sigma_{t-1}B \tag{5}$$

where $C$, $A$, and $B$ are $k \times k$ matrices but $C$ is upper triangular. An advantage of the BEKK model is that $\Sigma_t$ is positive definite if the diagonal elements of $C$ is positive. On the other hand, the model contains many parameters that do not represent directly the impact of $e_{t-1}$ or $\Sigma_{t-1}$ on the elements of $\Sigma_t$. In other words, it is hard to interpret the parameters of a BEKK model. Limited experience also shows that many parameter estimates of the BEKK model in Eq. (5) are statistically insignificant, implying that the model is overparameterized.

Using the standardization of Eq. (3), one can divide the multivariate volatility modeling into two steps. The first step is to specify models for elements of the $D_t$ matrix, and the second step is to model the correlation matrix $R_t$. Two such approaches have been proposed in the literature. In both cases, the elements $\sigma_{ii,t}$ are assumed to follow a univariate GARCH model. In other words, $\sigma_{ii,t}$ are based entirely on the $i$-the return series.

*Dynamic correlation model of Tse and Tsui.* In [8], the authors propose that (a) the individual volatility $\sigma_{ii,t}$ can assume any univariate GARCH models, and (b) the correlation matrix $R_t$ of Eq. (3) follows the model

$$R_t = (1 - \lambda_1 - \lambda_2)R + \lambda_1 \Psi_{t-1} + \lambda_2 R_{t-1} \tag{6}$$

where $\lambda_1$ and $\lambda_2$ are non-negative parameters satisfying $0 \leq \lambda_1 + \lambda_2 < 1$, $R$ is a $k \times k$ positive definite parameter matrix with $R_{ii} = 1$ and $\Psi_{t-1}$ is the $k \times k$ correlation matrix of some recent asset returns. For instance, if the most recent $m$ returns are used to define $\Psi_{t-1}$, then the $(i,j)$th element of $\Psi_{t-1}$ is given by

$$\psi_{ij,t-1} = \frac{\sum_{v=1}^{m} u_{i,t-v} u_{j,t-v}}{\sqrt{(\sum_{v=1}^{m} u_{i,t-v}^2)(\sum_{v=1}^{m} u_{j,t-v}^2)}},$$

where $u_{it} = e_{it}/\sqrt{\sigma_{ii,t}}$. If $m > k$, then $\Psi_{t-1}$ is positive definite almost surely. This in turn implies that $R_t$ is positive definite almost surely. We refer to this model as a $\text{DCC}_T(m)$ model. In practice, one can use the sample correlation matrix of the data to estimate $R$ in order to simplify the calculation. Indeed, this is the approach we shall take in this paper.

From the definition, the use of $\text{DCC}_T(m)$ model involves two steps. In the first step, univariate GARCH models are built for each return series. At step 2, the correlation matrix $R_t$ of Eq. (6) is estimated across all return series via the maximum likelihood method. An advantage of the $\text{DCC}_T(m)$ model is that the resulting correlation matrices are positive definite almost surely. In addition, the model is parsimonious in parameterization because the evolution of correlation matrices is governed by two parameters. On the other hand, strong limitation is imposed on the time evolution of the correlation matrices. In addition, it is hard to interpret the results of the two-step estimation. For instance, it is not clear what is the joint distribution of the innovation $e_t$ of the return series.

*Dynamic correlation model of Engle.* A similar correlation model is proposed by Engle [4]. Here the correlation matrix $R_t$ follows the model

$$R_t = W_t^{-1} Q_t W_t^{-1} \tag{7}$$

where $Q_t = [q_{ij,t}]$ is a positive-definite matrix, $W_t = \text{diag}\{\sqrt{q_{11,t}}, \ldots, \sqrt{q_{kk,t}}\}$ is a normalization matrix, and the elements of $Q_t$ are given by

$$Q_t = (1 - \alpha_1 - \alpha_2)\bar{Q} + \alpha_1 u_{t-1} u'_{t-1} + \alpha_2 Q_{t-1},$$

where $u_t$ is the standardized innovation vector with elements $u_{it} = e_{it}/\sqrt{\sigma_{ii,t}}$, $\bar{Q}$ is the sample covariance matrix of $u_t$, and $\alpha_1$ and $\alpha_2$ are non-negative scalar parameters satisfying $0 < \alpha_1 + \alpha_2 < 1$. We refer to this model as the $\text{DCC}_E$ model.

Compared with the $\text{DCC}_T(m)$ model, the $\text{DCC}_E$ model only uses the most recent standardized innovation to update the time-evolution of the correlation matrix. Since $u_{t-1}u'_{t-1}$ is singular for $k > 1$ and is, in general, not a correlation matrix, and the matrix $Q_t$ must be normalized in Eq. (7) to ensure that $R_t$ is indeed a correlation matrix. Because a single innovation is more variable than the correlation matrix of $m$ standardized innovations, the correlations of a $\text{DCC}_E$ model tend to be more variable than those of a $\text{DCC}_T(m)$ model.

To better understand the difference between $\text{DCC}_T(m)$ and $\text{DCC}_E$ models, consider the correlation $\rho_{12,t}$ of the first two returns in $r_t$. For $\text{DCC}_T(m)$ model,

$$\rho_{12,t} = (1 - \lambda_1 - \lambda_2)\rho_{12} + \lambda_2 \rho_{12,t-1} + \lambda_1 \frac{\sum_{v=1}^m u_{1,t-v} u_{2,t-v}}{\sqrt{(\sum_{v=1}^m u_{1,t-v}^2)(\sum_{v=1}^m u_{2,t-v}^2)}}.$$

On the other hand, for the $\text{DCC}_E$ model,

$$\rho_{12,t} = \frac{\alpha^* \bar{q}_{12} + \alpha_1 u_{1,t-1} u_{2,t-1} + \alpha_2 q_{12,t-1}}{\sqrt{(\alpha^* \bar{q}_{11} + \alpha_1 u_{1,t-1}^2 + \alpha_2 q_{11,t-1})(\alpha^* \bar{q}_{22} + \alpha_1 u_{2,t-1}^2 + \alpha_2 q_{22,t-1})}},$$

where $\alpha^* = 1 - \alpha_1 - \alpha_2$. The difference is clearly seen.

## 3. Proposed models

We start with the simple case in which the effects of positive and negative past returns on the volatility are symmetric. The case of asymmetric effects is given later.

### 3.1. Multivariate GARCH models

In this paper, we propose the following model

(8) $$r_t = \mu_t + e_t, \quad \mu_t = \phi_0 + \sum_{i=1}^p \phi_i r_{t-i}, \quad e_t = \Sigma^{1/2} \epsilon_t$$

where $p$ is a non-negative integer and $\{\epsilon_t\}$ is a sequence of independent and identically distributed multivariate Student-$t$ distribution with $v$ degrees of freedom. The probability density function of $\epsilon_t$ is

$$f(\epsilon) = \frac{\Gamma((v+k)/2)}{[\pi(v-2)]^{k/2} \Gamma(v/2)} [1 + (v-2)^{-1} \epsilon' \epsilon]^{-(v+k)/2}.$$

The variance of each component of $\epsilon_t$ is 1. The volatility matrix is standardized as Eq. (3) with elements of $D_t$ and the correlation matrix $R_t$ satisfying

(9) $$D_t^2 = \Lambda_0 + \Lambda_1 D_{t-1}^2 + \Lambda_2 G_{t-1}^2,$$
(10) $$R_t = (1 - \theta_1 - \theta_2)\bar{R} + \theta_1 \psi_{t-1} + \theta_2 R_{t-1},$$

where $G_t = \text{diag}\{e_{1t}, \ldots, e_{kt}\}$, $\Lambda_i = \text{diag}\{\ell_{11,i}, \ldots, \ell_{kk,i}\}$ are diagonal matrices such that $\ell_{ii,1} + \ell_{ii,2} < 1$ and $0 \leq \ell_{ii,j}$ for $i = 1, \ldots, k$ and $j = 1, 2$, $\bar{R}$ is the sample correlation matrix, $\theta_i$ are non-negative real numbers satisfying $\theta_1 + \theta_2 < 1$, and $\psi_{t-1}$ is the sample correlation matrix of the last $m$ innovations as defined in the $\text{DCC}_T(m)$ model of Eq. (6). We use $m = k + 2$ in empirical data analysis.

This model uses univariate GARCH(1,1) models for the conditional variance of components of $r_t$ and a combination of the correlation matrix equations of the $\text{DCC}_T(m)$ and $\text{DCC}_E$ models for the correlation. The order of GARCH models can be increased if necessary, but we use (1,1) for simplicity. In addition, $\Lambda_1$ and $\Lambda_2$ can be generalized to non-diagonal matrices. However, we shall keep the simple structure in Eq. (9) and (10) for ease in application and interpretation.

The proposed model differs from the $\text{DCC}_T(m)$ model in several ways. First, the proposed model uses a multivariate Student-$t$ distribution for innovation so that the degrees of freedom are the same for all asset returns. This simplifies the model interpretation at the expense of decreased flexibility. Second, the proposed model uses sample correlation matrix $\bar{R}$ to reduce the number of parameters. Third, the proposed model uses joint estimation whereas the $\text{DCC}_T(m)$ model performs separate estimations for variances and correlations.

### 3.2. Model with leverage effects

In financial applications, large positive and negative shocks to an asset have different impacts on the subsequent price movement. In volatility modeling, it is expected that a large negative shock would lead to increased volatility as a big drop in asset price is typically associated with bad news which, in turn, means higher risk for the investment. This phenomenon is referred to as the *leverage* effect in volatility modeling. The symmetry of GARCH model in Eq. (9) keeps the model simple, but fails to address the leverage effect. To overcome this shortcoming, we consider the modified model

$$(11) \qquad D_t^2 = \Lambda_0 + \Lambda_1 D_{t-1}^2 + \Lambda_2 G_{t-1}^2 + \Lambda_3 L_{t-1}^2,$$

where $\Lambda_i$ ($i = 0, 1, 2$) are defined as before, $\Lambda_3$ is a $k \times k$ diagonal matrix with non-negative diagonal elements, and $L_{t-1}$ is also a $k \times k$ diagonal matrix with diagonal elements

$$L_{ii,t-1} = \begin{cases} e_{i,t-1} & \text{if } e_{i,t-1} < 0, \\ 0 & \text{otherwise.} \end{cases}$$

In Eq. (11), we assume that $0 < \sum_{j=1}^{3} \ell_{ii,j} \leq 1$ for $i = 1, \ldots, k$. This is a sufficient condition for the existence of volatility.

From the definition, a positive shock $e_{i,t-1}$ affects the volatility via $\ell_{ii,2} e_{i,t-1}^2$. A negative shock, on the other hand, contributes $(\ell_{ii,2} + \ell_{ii,3}) e_{i,t-1}^2$ to the volatility. Checking the significance of $\ell_{ii,3}$ enables us to draw inference on the leverage effect.

### 4. Application

We illustrate the proposed model by considering some daily asset returns. First, we consider a four-dimensional process consisting of two equity returns and two exchange rate returns. Second, we consider a 10-dimensional equity returns. In both examples, we use $m = k + 2$ to estimate the local correlation matrices $\psi_{t-1}$ in Eq. (10).

**Example 1.** In this example, we consider the daily exchange rates between U.S. Dollar versus European Euro and Japanese Yen and the stock prices of IBM and Dell from January 1999 to December 2004. The exchange rates are the noon spot rate obtained from the Federal Reserve Bank of St. Louis and the stock returns are from the Center for Research in Security Prices (CRSP). We compute the simple returns of the exchange rates and remove returns for those days when one of the markets was not open. This results in a four-dimensional return series with 1496 observations. The return vector is $r_t = (r_{1t}, r_{2t}, r_{3t}, r_{4t})'$ with $r_{1t}$ and $r_{2t}$ being the returns of Euro and Yen exchange rate, respectively, and $r_{3t}$ and $r_{4t}$ are the returns of IBM and Dell stock, respectively. All returns are in percentages. Figure 1 shows the time plot of the return series. From the plot, equity returns have higher variability than the exchange rate returns, and the variability of equity returns appears to be decreasing in recent years. Table 1 provides some descriptive statistics of the return series. As expected, the means of the return are essentially zero and all four series have heavy tails with positive excess kurtosis.

The equity returns have some serial correlations, but the magnitude is small. If multivariate Ljung-Box statistics are used, we have $Q(3) = 59.12$ with p-value 0.13 and $Q(5) = 106.44$ with p-value 0.03. For simplicity, we use the sample mean as the mean equation and apply the proposed multivariate volatility model to the mean-corrected data. In estimation, we start with a general model, but add some equality constraints as some estimates appear to be close to each other. The results are given in Table 2 along with the value of likelihood function evaluated at the estimates.

For each estimated multivariate volatility model in Table 2, we compute the

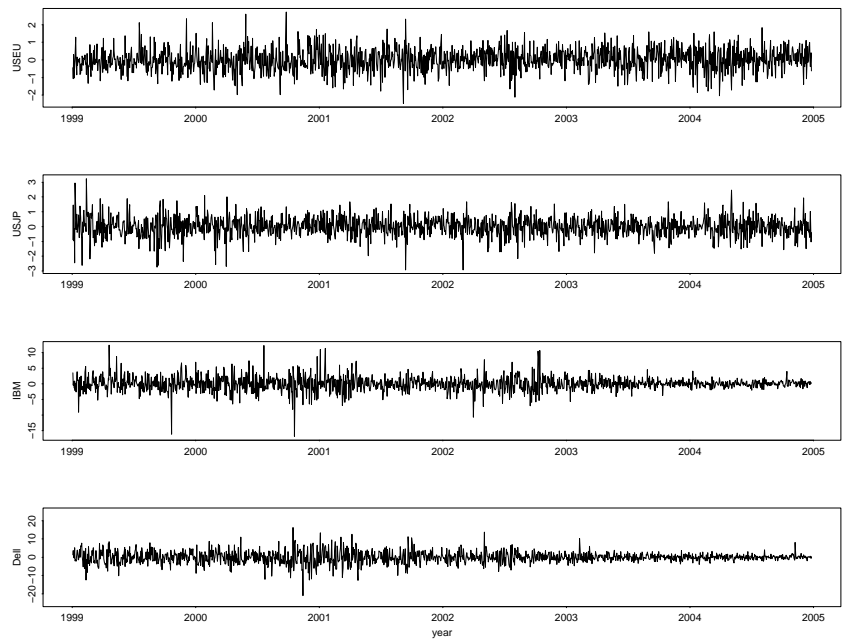

FIG 1. *Time plots of daily return series from January 1999 to December 2004:* (a) *Dollar-Euro exchange rate,* (b) *Dollar-Yen exchange rate,* (c) *IBM stock, and* (d) *Dell stock.*

TABLE 1
*Descriptive statistics of daily returns of Example 1. The returns are in percentages, and the sample period is from January 1999 to December 2004 for 1496 observations*

| Asset | USEU | JPUS | IBM | DELL |
|---|---|---|---|---|
| Mean | 0.0091 | −0.0059 | 0.0066 | 0.0028 |
| Standard error | 0.6469 | 0.6626 | 5.4280 | 10.1954 |
| Skewness | 0.0342 | −0.1674 | −0.0530 | −0.0383 |
| Excess kurtosis | 2.7090 | 2.0332 | 6.2164 | 3.3054 |
| Box-Ljung $Q(12)$ | 12.5 | 6.4 | 24.1 | 24.1 |

TABLE 2
*Estimation results of multivariate volatility models for Example 1 where $L_{\max}$ denotes the value of likelihood function evaluated at the estimates, $v$ is the degrees of freedom of the multivariate Student-t distribution, and the numbers in parentheses are asymptotic standard errors*

| $\Lambda_0$ | $\Lambda_1$ | $\Lambda_2$ | $(v, \theta_1, \theta_2)'$ |
|---|---|---|---|
| (a) Full model estimation with $L_{\max} = -9175.80$ | | | |
| 0.0041(0.0033) | 0.9701(0.0114) | 0.0214(0.0075) | 7.8729(0.4693) |
| 0.0088(0.0038) | 0.9515(0.0126) | 0.0281(0.0084) | 0.9808(0.0029) |
| 0.0071(0.0053) | 0.9636(0.0092) | 0.0326(0.0087) | 0.0137(0.0025) |
| 0.0150(0.0136) | 0.9531(0.0155) | 0.0461(0.0164) | |
| (b) Restricted model with $L_{\max} = -9176.62$ | | | |
| 0.0066(0.0028) | 0.9606(0.0068) | 0.0255(0.0068) | 7.8772(0.7144) |
| 0.0066(0.0023) | | 0.0240(0.0059) | 0.9809(0.0042) |
| 0.0080(0.0052) | | 0.0355(0.0068) | 0.0137(0.0025) |
| 0.0108(0.0086) | | 0.0385(0.0073) | |
| (c) Final restricted model with $L_{\max} = -9177.44$ | | | |
| 0.0067(0.0021) | 0.9603(0.0063) | 0.0248(0.0048) | 7.9180(0.6952) |
| 0.0067(0.0021) | | 0.0248(0.0048) | 0.9809(0.0042) |
| 0.0061(0.0044) | | 0.0372(0.0061) | 0.0137(0.0028) |
| 0.0148(0.0084) | | 0.0372(0.0061) | |
| (d) Model with leverage effects, $L_{\max} = -9169.04$ | | | |
| 0.0064(0.0027) | 0.9600(0.0065) | 0.0254(0.0063) | 8.4527(0.7556) |
| 0.0066(0.0023) | | 0.0236(0.0054) | 0.9810(0.0044) |
| 0.0128(0.0055) | | 0.0241(0.0056) | 0.0132(0.0027) |
| 0.0210(0.0099) | | 0.0286(0.0062) | |

standardized residuals as

$$\hat{\epsilon}_t = \widehat{\Sigma}_t^{-1/2} e_t,$$

where $\widehat{\Sigma}_t^{1/2}$ is the symmetric square-root matrix of the estimated volatility matrix $\widehat{\Sigma}_t$. We apply the multivariate Ljung-Box statistics to the standardized residuals $\hat{\epsilon}_t$ and its squared process $\hat{\epsilon}_t^2$ of a fitted model to check model adequacy. For the full model in Table 2(a), we have $Q(10) = 167.79(0.32)$ and $Q(10) = 110.19(1.00)$ for $\hat{\epsilon}_t$ and $\hat{\epsilon}_t^2$, respectively, where the number in parentheses denotes p-value. Clearly, the model adequately describes the first two moments of the return series. For the model in Table 2(b), we have $Q(10) = 168.59(0.31)$ and $Q(10) = 109.93(1.00)$. For the final restricted model in Table 2(c), we obtain $Q(10) = 168.50(0.31)$ and $Q(10) = 111.75(1.00)$. Again, the restricted models are capable of describing the mean and volatility of the return series.

From Table 2, we make the following observations. First, using the likelihood ratio test, we cannot reject the final restricted model compared with the full model. This results in a very parsimonious model consisting of only 9 parameters for the time-varying correlations of the four-dimensional return series. Second, for the two stock return series, the constant terms in $\Lambda_0$ are not significantly different from zero and the sum of GARCH parameters is $0.0372 + 0.9603 = 0.9975$, which is very close to unity. Consequently, the volatility series of the two equity returns exhibit IGARCH behavior. On the other hand, the volatility series of the two exchange rate returns appear to have a non-zero constant term and high persistence in GARCH parameters. Third, to better understand the efficacy of the proposed model, we compare the results of the final restricted model with those of rolling estimates. The rolling estimates of covariance matrix are obtained using a moving window of size 69, which is the approximate number of trading days in a quarter. Figure 2 shows the time plot of estimated volatility. The solid line is the volatility obtained by the proposed model and the dashed line is for volatility of the rolling estimation. The overall pattern seems similar, but, as expected, the rolling estimates respond slower than the proposed model to large innovations. This is shown by the faster rise and decay of the volatility obtained by the proposed model. Figure 3 shows the time-varying correlations of the four asset returns. The solid line denotes correlations obtained by the final restricted model of Table 2 whereas the dashed line is for rolling estimation. The correlations of the proposed model seem to be smoother.

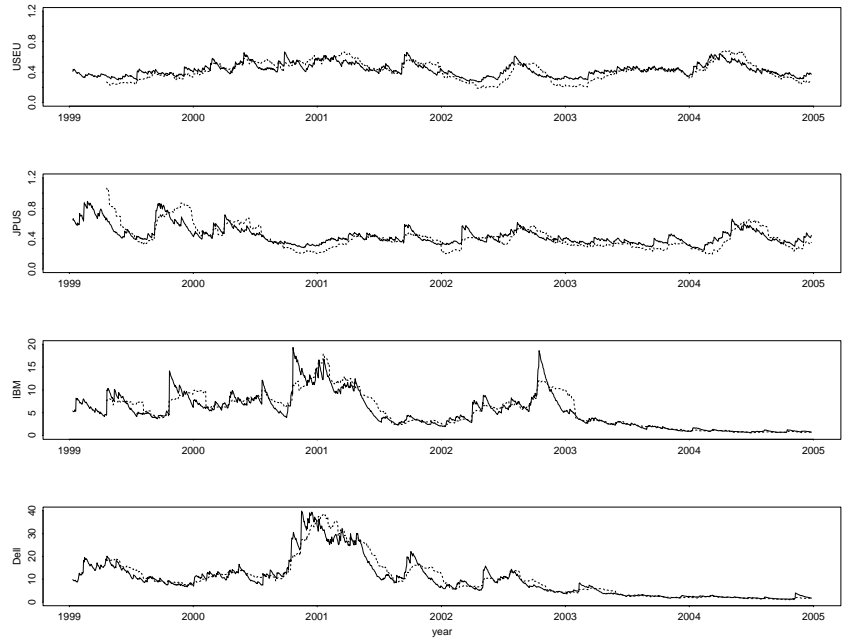

Fig 2. *Time plots of estimated volatility series of four asset returns. The solid line is from the proposed model and the dashed line is from a rolling estimation with window size* 69: (a) *Dollar-Euro exchange rate,* (b) *Dollar-Yen exchange rate,* (c) *IBM stock, and* (d) *Dell stock.*

Table 2(d) gives the results of a fitted integrated GARCH-type with leverage effects. The leverage effects are statistically significant for equity returns only and are in the form of an IGARCH model. Specifically, the $\Lambda_3$ matrix of the correlation equation in Eq. (11) is

$$\Lambda_3 = \text{diag}\{0, 0, (1 - 0.96 - 0.0241), (1 - 0.96 - 0.0286)\} = \text{diag}\{0, 0, 0.0159, 0.0114\}.$$

Although the magnitudes of the leverage parameters are small, but they are statistically significant. This is shown by the likelihood ratio test. Specifically, compared the fitted models in Table 2(b) and (d), the likelihood ratio statistic is 15.16, which has a p-value 0.0005 based on the chi-squared distribution with 2 degrees of freedom.

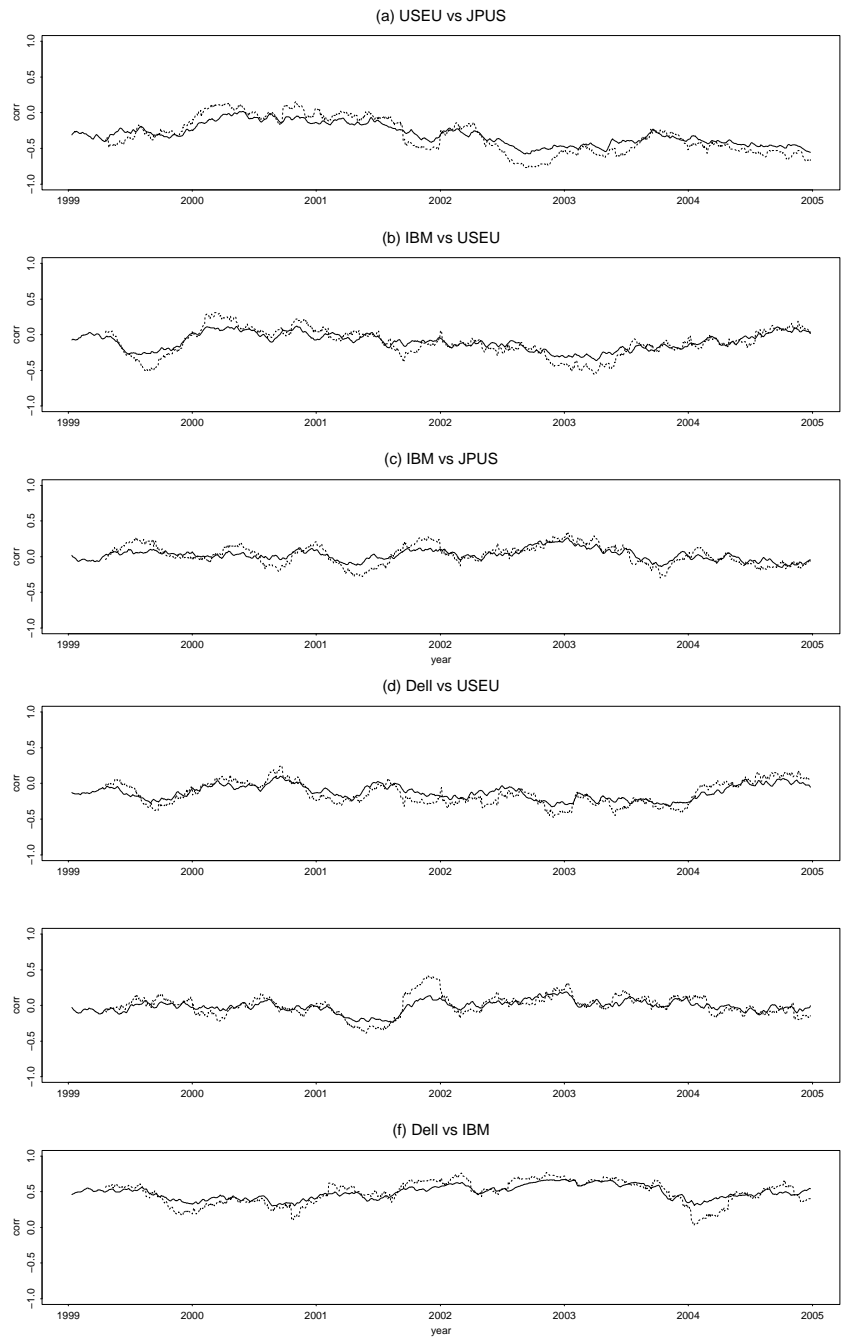

FIG 3. *Time plots of time-varying correlations between the percentage simple returns of four assets from January 1999 to December 2004. The solid line is from the proposed model whereas the dashed line is from a rolling estimation with window size* 69.

**Example 2.** In this example, we consider daily simple returns, in percentages, of the S&P 500 index and nine individual stocks from January 1990 to December 2004 for 3784 observations. Thus, we have a 10-dimensional return series. The ten assets are given in Table 3 along with some descriptive statistics. All asset returns have positive excess kurtosis, indicating heavy tails. Except for the stock of General Electrics, return minimums exceed the maximums in modulus, suggesting asymmetry in price changes due to good and bad news.

Sincere there are some minor serial and cross correlations in the 10-dimensional returns, we fit a vector autoregressive model of order 3, i.e. VAR(3), to the data to remove the dynamic dependence and employ the resulting residual series in volatility modeling. See Eq. (8).

We have applied the proposed volatility model in Eqs. (9)- (10) to the residual series of the VAR(3) model. But our subsequently analysis shows that the model with leverage effects in Eq. (11) is preferred based on the likelihood ratio test. Therefore, for simplicity in presentation, we shall only report the results with leverage effects.

Employing the volatility model in Eq. (11) with the correlations in Eq. (10), we found that for the returns of IBM, DELL, GE, and GM stocks the leverage effects follow integrated GARCH models. Consequently, for these four stock returns the leverage parameters are given by

$$\Lambda_{ii,3} = 1 - \Lambda_{ii,1} - \Lambda_{ii,2},$$

where $\Lambda_{ii,j}$ is the $i$th diagonal element of the matrix $\Lambda_j$, $j = 1, 2, 3$. Table 4 shows the parameter estimates of the 10-dimensional volatility model.

For model checking, we use a bootstrap method to generate the critical values of multivariate Ljung-Box statistics for the standardized residuals and their squared series. Specifically, we generate 10,000 realizations each with 3781 observations from the standardized residuals of the fitted model. The bootstrap samples are drawn with replacement. For each realization, we compute the Ljung-Box statistics $Q(5)$, $Q(10)$, and $Q(15)$ of the series and its squared series. Table 5 gives some critical values of the Ljung-Box statistics. For the fitted model, we have $Q(10) = 836.12$ and $Q(15) = 1368.71$ for the standardized residuals and $Q(10) = 1424.64$ and $Q(15) = 1923.75$ for the squared series of standardized residuals. Compared with the critical values in Table 5, the Ljung-Box statistics are not significant at the 1% level. Thus, the fitted model is adequate in modeling the volatility of the 10-dimensional return series. We also applied the AIC criterion to the squared series of standardized residuals. The criterion selected a VAR(0) model, confirming that the fitted multivariate volatility model is adequate.

From the fitted model, we make the following observations. First, except for two marginal cases, all estimates of leverage parameters are statistically significant at the 5% level based on their asymptotic standard errors. The two marginally significant leverage parameters are for BA amd PFE stocks and their $t$-ratios are 1.65 and 1.92, respectively. Thus, as expected, the leverage effects are positive and most of them are significant. Second, all parameters of the volatility equation are significant. Thus, the model does not contain unnecessary parameters. Third, the model contains 30 parameters. This is very parsimonious for a 10-dimensional return series. Fourth, the correlations evolve slowly with high persistence parameter 0.9864.

TABLE 3
*Descriptive statistics of asset returns used in Example 2. Except for the S&P index, tick symbol is used to denote the company. Returns are in percentages*

| Asset | Mean | St.Error | Skewness | Ex.Kurt. | Minimum | Maximum |
|-------|------|----------|----------|----------|---------|---------|
| S&P   | 0.038 | 1.03 | −0.018 | 3.58 | −6.87 | 5.73 |
| IBM   | 0.066 | 2.03 | 0.294  | 6.32 | −15.54 | 13.16 |
| INTC  | 0.122 | 2.82 | −0.122 | 4.17 | −22.02 | 20.12 |
| DELL  | 0.236 | 3.49 | −0.012 | 3.45 | −25.44 | 20.76 |
| GE    | 0.074 | 1.70 | 0.176  | 3.80 | −10.67 | 12.46 |
| BA    | 0.052 | 1.98 | −0.282 | 6.08 | −17.62 | 11.63 |
| GM    | 0.039 | 2.01 | 0.111  | 1.98 | −13.53 | 10.34 |
| JNJ   | 0.076 | 1.59 | −0.139 | 4.32 | −15.85 | 8.21 |
| MRK   | 0.051 | 1.80 | −0.861 | 14.91 | −26.78 | 9.60 |
| PFE   | 0.084 | 1.91 | −0.068 | 1.94 | −11.15 | 9.71 |

TABLE 4
*Parameter estimates of the proposed volatility model with leverage effects for the 10 asset returns of Example 2. For leverage effects, those estimates without standard errors denote IGARCH constraints*

| $\Lambda_1$ | | | | $\lambda I$ | | | | | | |
|---|---|---|---|---|---|---|---|---|---|---|
| Estimate | | | | 0.9658 | | | | | | |
| Std.Err | | | | 0.0024 | | | | | | |
| $\Lambda_2$ | | | | Diagonal matrix with elements | | | | | | |
| Estimate | .0154 | .0174 | .0168 | .0298 | .0191 | .0206 | .0187 | .0110 | .0128 | .0192 |
| Std.Err | .0031 | .0026 | .0038 | .0030 | .0029 | .0041 | .0037 | .0038 | .0028 | .0037 |
| $\Lambda_0$ | | | | Diagonal matrix with elements | | | | | | |
| Estimate | .0077 | .0211 | .0763 | .0170 | .0185 | .0279 | .0342 | .0281 | .0369 | .0309 |
| Std.Err | .0010 | .0042 | .0121 | .0067 | .0031 | .0054 | .0074 | .0048 | .0061 | .0058 |
| $\Lambda_3$ | | | | Diagonal matrix with elements | | | | | | |
| Estimate | .0178 | .0168 | .0126 | .0044 | .0152 | .0107 | .0155 | .0210 | .0143 | .0115 |
| Std.Err | .0049 | | .0059 | | | .0065 | | .0064 | .0059 | .0060 |
| Parameter | $v$ | $\theta_2$ | $\theta_1$ | | | | | | | |
| Estimate | 9.54 | .9864 | .0070 | | | | | | | |
| Std.Err | .417 | .0016 | .0006 | | | | | | | |

TABLE 5
*Critical values of Ljung-Box statistics for 10-dimensional standardized residual series. The values are obtained by a bootstrap method with 10,000 iterations. The sample size of the series is 3781*

| | Standardized residuals | | | Squared standardized residuals | | |
|---|---|---|---|---|---|---|
| Statistics | 1% | 5% | 10% | 1% | 5% | 10% |
| $Q(5)$ | 576.92 | 553.68 | 541.33 | 915.89 | 696.82 | 617.74 |
| $Q(10)$ | 1109.05 | 1075.25 | 1057.94 | 1150.31 | 1281.03 | 1170.12 |
| $Q(15)$ | 1633.31 | 1591.61 | 1571.17 | 2125.65 | 1837.28 | 1713.50 |

Fifth, the estimate of the degrees of freedom for multivariate Student-$t$ innovation is 9.54, confirming that the returns have heavy tails.

**Remark.** In this paper, we use a MATLAB program to estimate the proposed multivariate volatility models. The negative log likelihood function is minimized with inequality parameter constraints using the function **fmincon**. Limited experience shows that the results are not sensitive to the initial values, but initial values that are far away from the final estimates do require many more iterations. The estimation, however, can become difficult if some parameters are approching the boundary of the parameter space. For instance, if there is no leverage effect, then the hessian matrix can be unstable when the leverage parameter is included in the estimation.

## 5. Extensions and some alternative approaches

In this paper, we consider a simple approach to model multivariate volatilities of asset returns. Unlike other methods available in the literature, the proposed approach estimates the conditional variances and correlations jointly and the resulting volatility matrices are positive definite. The proposed model can handle the leverage effects and is parsimonious. We demonstrated the efficacy of the proposed model by analyzing a 4-dimensional and a 10-dimensional asset return series. The results are encouraging. We also used a bootstrap method to obtain finite-sample critical values for the multivariate Ljung-Box statistics for testing serial and cross correlations of a vector series.

There are possible extensions of the proposed model. For example, Eq. (10) requires that all correlations have the same persistence parameter $\theta_2$. This restriction can be relaxed by letting $\theta_1$ and $\theta_2$ be diagonal matrices of positive real numbers. The model would become

$$R_t = (I - \theta_1^2 - \theta_2^2)\bar{R} + \theta_1 \psi_{t-1} \theta_1 + \theta_2 R_{t-1} \theta_2.$$

Under this model, the $i$th asset return contributes $\theta_{ii,2}$ to the persistence of correlations. In addition, one can have equality constraints among diagonal elements of each $\theta_i$ matrix to keep the model parsimonious.

Some alternative approaches have been considered in the literature to overcome the curse of dimensionality in multivariate volatility modeling. Palandri [7] uses a sequential Cholesky decomposition to build a multivariate volatility of 69 stock returns. The independent component models have also been used to simplify the modeling procedure, e.g., see [6].